\documentclass[11pt]{article}%
\usepackage{amsmath}
\usepackage{graphicx}%
\usepackage{amsfonts}%
\usepackage{amssymb}
\usepackage{makecell}
\usepackage{fancyhdr}
\newcommand{\eqnb}{\begin{equation}}
\newcommand{\eqne}{\end{equation}}

\newtheorem{The}{Theorem}

\newtheorem{Rem}{Remark}
\allowdisplaybreaks

\begin{document}

\title{A Nonlinear Solution to Closed Queueing Networks for Bike Sharing Systems with Markovian
Arrival Processes and under an Irreducible Path Graph}
\author{Quan-Lin Li, Rui-Na Fan and Zhi-Yong Qian\\School of Economics and Management Sciences \\Yanshan University, Qinhuangdao 066004, P.R. China}
\maketitle
\pagestyle{plain} 
\begin{abstract}
As a favorite urban public transport mode, the bike sharing system is a large-scale and complicated system,
and there exists a key requirement that a user and a bike should be matched sufficiently in time. Such matched behavior
makes analysis of the bike sharing systems more difficult and challenging. To design a better bike sharing system, it is a key to analyze and compute the probabilities of the problematic (i.e., full or empty) stations. In fact, such a computation is established for some fairly complex stochastic systems. To do this, this paper considers a more general large-scale bike sharing system from two
important views: (a) Bikes move in an irreducible path graph, which is related to geographical structure of the bike sharing system; and (b) Markovian arrival processes (MAPs) are applied to describe the non-Poisson and burst behavior of
bike-user (abbreviated as user) arrivals, while the burstiness demonstrates that the user arrivals are time-inhomogeneous and space-heterogeneous in practice. For such a complicated bike sharing system, this paper establishes a multiclass closed queueing network by means of
some virtual ideas, for example, bikes are abstracted as virtual customers; stations and roads are regarded as virtual nodes. Thus user arrivals
are related to service times at station nodes; and users riding bikes on roads are viewed as service times at road nodes.
Further, to deal with this multiclass closed queueing network, we provide a detailed observation practically on physical behavior of the bike sharing system in order to establish the routing matrix, which gives a nonlinear solution to compute the relative arrival rates in terms of the product-form solution to the steady-state probabilities of joint queue lengths at the virtual nodes. Based on this, we can compute the steady-state probability of problematic stations, and also deal with other interesting performance measures of the bike sharing system. We hope that the methodology and results of this paper can be applicable in the study of more general bike sharing systems through multiclass closed queueing networks.
\vskip                                                                 0.5cm


\end{abstract}

\section{Introduction}
In this paper, we propose a more general bike sharing system with Markovian arrival processes and under an irreducible path graph.
Note that the bike sharing system always has some practically important factors, for example, time-inhomogeneity, geographical heterogeneity, and arrival burstiness. To analyze such a bike sharing system, we establish a multiclass closed queueing network by means of virtual customers, virtual nodes and virtual service times. Further, when studying this multiclass closed queueing network, we set up a routing matrix which gives a nonlinear solution to compute the relative arrival rates, and provide the product-form solution to the steady-state probabilities of joint queue
lengths at the virtual nodes. Based on this, we can compute the steady-state probability of problematic stations, and also deal with
other interesting performance measures of the bike sharing system.

During the last decades bike sharing systems have emerged as a public transport
mode devoted to short trip in more than 600 major cities around the world.
Bike sharing systems are regarded as a promising way to jointly reduce, such as, traffic and parking
congestion, traffic noise, air pollution and greenhouse effect. Several excellent overviews and useful remarks were
given by DeMaio \cite{DeM:2009}, Meddin and DeMaio \cite{Med:2012}, Shu \emph{et al.} \cite{Shu:2013}, Labadi \emph{et al}. \cite{Lab:2015} and
Fishman \emph{et al}. \cite{Fis:2013}.

Few papers applied queueing theory and Markov processes to the study of bike sharing systems. On this research line, it is a key to compute the probability of problematic stations. However, so far there still exist some basic difficulties and challenges for computing the probability of problematic stations because computation of the steady-state probability, in the bike sharing system, needs to apply the theory of complicated or high-dimensional Markov processes. For this, readers may refer to recent literatures which are classified
and listed as follows. \textbf{(a) Simple queues:} Leurent \cite{Leu:2012} used the M/M/1/C queue to
study a vehicle-sharing system, and also analyzed performance measures
of this system. Schuijbroek \emph{et al}. \cite{Sch:2013} evaluated the service level by means of the transient
distribution of the M/M/1/C queue, and the service level was used to
establish some optimal models to discuss vehicle routing. Raviv \emph{et al}. \cite{Rav:2013} and Raviv and Kolka \cite{Rav:2013a}
employed the transient distribution of the time-inhomogeneous M$\left(
t\right)$/M$\left(  t\right)$/1/C queue to compute the expected number of
bike shortages at each station. \textbf{(b) Closed queueing networks:} Adelman
\cite{Ade:2007} applied a closed queueing network to propose an internal
pricing mechanism for managing a fleet of service units, and also used a
nonlinear flow model to discuss the price-based policy for establishing the vehicle
redistribution. George and Xia \cite{Geo:2011} used the closed queueing
networks to study the vehicle rental systems, and determined the optimal
number of parking spaces for each rental location. Li \emph{et al}. \cite{Li:2016a}
proposed a unified framework for analyzing the closed queueing networks in
the study of bike sharing systems. \textbf{(c) Mean-field method}. Fricker \emph{et al}.
\cite{Fri:2012} considered a space-inhomogeneous bike-sharing system with
multiple clusters, and expressed the minimal proportion of problematic
stations. Fricker and Gast \cite{Fri:2016} provided a detailed analysis for a
space-homogeneous bike-sharing system in terms of the M/M/1/K queue as well
as some simple mean-field models, and crucially, they derived the closed-form
solution to find the minimal proportion of problematic stations. Fricker and Tibi
\cite{Fri:2017} studied the central limit and local limit theorems for the
independent (non-identically distributed) random variables, which provide support on
analysis of a generalized Jackson network with product-form solution. Further,
they used the limit theorems to give an outline of stationary asymptotic
analysis for the locally space-homogeneous bike-sharing systems. Li \emph{et al}. \cite{Li:2016b} provided a complete picture on how to jointly use the
mean-field theory, the time-inhomogeneous queues and the nonlinear birth-death
processes to analyze performance measures of the bike-sharing systems. Li and Fan \cite{Li:2016} discussed the
bike sharing system under an Markovian environment by means of the mean-field computation,
the time-inhomogeneous queues and the nonlinear Markov processes. \textbf{(d) Markov decision processes}. To
discuss the bike-sharing systems, Waserhole and Jost \cite{Was:2012, Was:2013, Was:2016} and Waserhole \emph{et al}. \cite{Was:2013a} used
the simplified closed queuing networks to establish the Markov decision
models, and computed the optimal policy by means of the fluid approximation
which overcame the state space explosion of multi-dimensional Markov decision
processes.

There has been much key research on closed queueing networks. Readers may refer to, such as, three excellent books by Kelly \cite{Kel:1979, Kel:2011} and Serfozo \cite{Ser:1999}; multiclass customers by Baskett \emph{et al}. \cite{Bas:1975}, multiple closed chains by Reiser and Kobayashi \cite{Rei:1975}, computational algorithms by Bruell and Balbo \cite{Bru:1980}, mean-value computation by Reiser \cite{Rei:1981}, sojourn time by Kelly and Pollett \cite{Kel:1983}, survey for blocks by Onvural \cite{Onv:1990}, and batch service by Henderson \emph{et al}. \cite{Hen:1990}.

Markovian arrival process (MAP) is a useful mathematical model for describing bursty traffic in, for example, communication networks,
manufacturing systems, transportation networks and so forth. Readers may refer to recent publications for more details, among which are
Ramaswami \cite{Ram:1980}, Chapter 5 in Neuts \cite{Neu:1989}, Lucantoni \cite{Luc:1991}, Neuts \cite{Neu:1995},
Chakravarthy \cite{Cha:2001} and Li \cite{Li:2010}.

\textbf{Contributions of this paper:} The main contributions of this paper are twofold: The first contribution is to propose a more general bike sharing system with Markovian arrival processes and under an irreducible path graph. Note that Markovian arrival processes, as well as the irreducible path graph indicate that burst arrival behavior and geographical structure of the bike sharing system are more general and practical. Specifically, the burstiness is to well express that the user arrivals are time-inhomogeneous and space-heterogeneous in practice. For such a bike sharing system, this paper establishes a multiclass closed queueing network by means of virtual customers, virtual nodes and virtual service times. The second contribution is to deal with such a multiclass closed queueing network with virtual customers, virtual nodes and virtual service times, and to establish a routing matrix which gives a nonlinear solution to compute the relative arrival rates in terms of the product-form solution to the steady-state probabilities of joint queue lengths at the virtual nodes. By using the product-form solution, this paper computes the steady-state probability of problematic stations, and also deals with other interesting performance measures of the bike sharing system. Therefore, the methodology and results of this paper can be applicable in the study of more general bike sharing systems by means of multiclass closed queueing networks.

\textbf{Organization of this paper:} The remainder of this paper is organized as follows. In Section 2, we describe
a large-scale bike sharing system with Markovian arrival processes and under an irreducible path graph. In Section 3, we abstract the bike sharing system as a multiclass closed queueing network with virtual customers, virtual nodes and virtual service times. Further, we establish the routing matrix, and compute the relative arrival rate in each node, where three examples are given to express and compute the routing matrix and the relative arrival rate. In Section 4, we give a product-form solution to the steady-state probabilities of joint queue lengths at the virtual nodes, and provide
a nonlinear solution to determine the $N$ undetermined constants which are related to the probability of problematic stations. Moreover, we compute the steady-state probability of problematic stations, and also analyze other performance measures of the bike sharing system. Finally, some concluding remarks are given in Section 5.

\section{Model Description}

In this section, we describe a more general large-scale bike sharing system,
where arrivals of bike users are non-Poisson and are characterized as Markovian arrival
processes (MAPs), and users riding bikes travel in an irreducible path graph which is constituted by $N$ different stations and some different directed roads.

In a large-scale bike sharing system, a user arrives at a station, rents a
bike, and uses it for a while; then he returns the bike to another station,
and immediately leaves this system. Based on this, we describe a more general large-scale space-heterogeneous bike sharing system, and introduce operational mechanism, system
parameters and basic notation as follows:

\textbf{(1) Stations: }We assume that there are $N$ different stations in the
bike sharing system. The $N$\ stations may be different due to their
geographical location and surrounding environment. We assume that
every station has $C$ bikes and $K$ parking positions at the initial time
$t=0$, where $1\leq C<K<\infty$, and $NC\geq K$. Note that such a condition $NC\geq K$ is to make at
least a full station.

\textbf{(2) Roads: }Let Road $i\rightarrow j$ be a road relating Station $i$ to Station $j$. Note that Road $i\rightarrow j$ and Road $j\rightarrow i$ may be different. To express all the roads beginning from Station $i$ for $1\leq i\leq N$, we write
\[
R\left(  i\right)  =\left\{  \text{Road }i\rightarrow j:j\neq i,1\leq j\leq
N\right\}  .
\]
Similarly, to express all roads over at Station $j$ for $1\leq j\leq N$, we write
\begin{equation*}
\overline{R}\left( j\right) =\left\{ \text{Road }i\rightarrow j:i\neq
j,1\leq i\leq N\right\}.
\end{equation*}
It is easy to see that there are at most $N-1$ different directed roads in the set $R\left(
i\right)  $ or $\overline{R}\left( j\right)$. We denote by $\left|  R\left(  i\right)  \right|  $ the number of
elements or roads in the set $R\left(  i\right)  $. Thus $\left|  R\left(  i\right)
\right|  \leq N-1$ for $1\leq i\leq N$ and $\sum\nolimits_{i=1}^{N}\left|
R\left(  i\right)  \right|  \leq N\left(  N-1\right)  $.

To express all the stations in the near downlink of Station $i$, we write%
\[
\Theta_{i}=\left\{  j:\text{Road }i\rightarrow j\in R\left(  i\right)
\right\}  .
\]
Similarly, the set of all stations in the near uplink of Station $i$ is written as
\begin{equation*}
\Delta _{i}=\left\{ j:\text{Road }j\rightarrow i\in \overline{R}\left(
i\right) \right\} .
\end{equation*}

\textbf{(3) An irreducible path graph:} To express the bike moving paths, it is easy to observe that the bikes dynamically move among the stations and among the roads. To record the bike dynamic positions, it is better to introduce two classes
of virtual nodes: (a) station nodes; and (b) road nodes. The set of all
the virtual nodes of the bike sharing system is given by%
\begin{equation*}
\Theta =\left\{ \text{Station }i:1\leq i\leq N\right\} \cup \left\{ \underset%
{i=1}{\overset{N}{\cup }}R\left( i\right) \right\}.
\end{equation*}
In this bike sharing system, it is easy to calculate that there are $N+\cup_{i=1}^{N}\left|  R\left(
i\right)  \right|  $ virtual nodes.

If Station $i$ has a near downstream Road $i\rightarrow j$, then we call that
Node $i$ (i.e. Station $i$) can be accessible to Node $i\rightarrow j$
(i.e. Road $i\rightarrow j$), denoted as Node $i\Longrightarrow$ Node
$i\rightarrow j$; otherwise Node $i$ can not be accessible to Node
$i\rightarrow j$. If Station $j$ has a near upstream Road $i\rightarrow j$, then
we call that Node $i\rightarrow j$ can be accessible to Node $j$, denoted
as Node $i\rightarrow j\Longrightarrow$ Node $j$; otherwise Node
$i\rightarrow j$ can not be accessible to Node $j$.

If there exist some virtual nodes $n_{1},n_{2},\ldots,n_{r}$ in the set
$\Theta$ such that%
\[
\text{Node }n_{1}\Longrightarrow\text{Node }n_{2}\Longrightarrow
\cdots\Longrightarrow\text{Node }n_{r},
\]
then we call that there is an accessible path formed by the virtual nodes
$n_{1},n_{2},$ $\ldots,n_{r}$.

If for any two virtual nodes $m_{a}$ and $m_{b}$
in the set $\Theta$, there always exist some virtual nodes $n_{1}%
,n_{2},\ldots,n_{r}$ in the set $\Theta$ such that%
\[
\text{Node }m_{a}\Longrightarrow\text{Node }n_{1}\Longrightarrow
\text{Node }n_{2}\Longrightarrow\cdots\Longrightarrow\text{Node }%
n_{r}\Longrightarrow\text{Node }m_{b},
\]
then we call that the path graph of the bike sharing system is irreducible.

\emph{In this paper, we assume that the bike sharing system exists an irreducible path graph}. In this case,
we call that the bike sharing system is path irreducible. Note that this irreducibility is guaranteed through setting up an appropriate road construction with
$R\left(  i\right)  $ for $1\leq i\leq N$. In general, such a road construction is not unique in order to
guarantee the irreducible path graph.

\textbf{(4) Markovian arrival processes: }Arrivals of outside bike users
at Station $i$ are a Markovian arrival process (MAP) of irreducible matrix descriptor $\left(
\mathbf{C}_{i},\mathbf{D}_{i}\right)  $ of size $m$, denoted as MAP$\left(
\mathbf{C}_{i},\mathbf{D}_{i}\right)  $, where%

\[
\mathbf{C}_{i}=\left(
\begin{array}
[c]{cccc}%
c_{1,1}^{\left(  i\right)  } & c_{1,2}^{\left(  i\right)  } & \cdots &
c_{1,m}^{\left(  i\right)  }\\
c_{2,1}^{\left(  i\right)  } & c_{2,2}^{\left(  i\right)  } & \cdots &
c_{2,m}^{\left(  i\right)  }\\
\vdots & \vdots & \ddots & \vdots\\
c_{m,1}^{\left(  i\right)  } & c_{m,2}^{\left(  i\right)  } & \cdots &
c_{m,m}^{\left(  i\right)  }%
\end{array}
\right)
\]
and%

\[
\mathbf{D}_{i}=\left(
\begin{array}
[c]{cccc}%
d_{1,1}^{\left(  i\right)  } & d_{1,2}^{\left(  i\right)  } & \cdots &
d_{1,m}^{\left(  i\right)  }\\
d_{2,1}^{\left(  i\right)  } & d_{2,2}^{\left(  i\right)  } & \cdots &
d_{2,m}^{\left(  i\right)  }\\
\vdots & \vdots & \ddots & \vdots\\
d_{m,1}^{\left(  i\right)  } & d_{m,2}^{\left(  i\right)  } & \cdots &
d_{m,m}^{\left(  i\right)  }%
\end{array}
\right)  .
\]
Let $c_{k,l}^{\left(i\right)  }\geq0$ with $l\neq k$, $d_{r,s}^{\left(i\right)  }\geq0$, $c_{k,k}^{\left( i\right) }=-\left( \sum\limits_{l\neq k}^{m}c_{k,l}^{\left(
i\right) }+\sum\limits_{r=1}^{m}d_{k,r}^{\left( i\right) }\right) $, and hence $\left(
\mathbf{C}_{i}+\mathbf{D}_{i}\right)  e=0$. We assume that Markov chain
$\mathbf{C}_{i}+\mathbf{D}_{i}$ is irreducible, finite-state and aperiodic, hence it is positive-recurrent due to the finite state space. Further, in the
Markov chain $\mathbf{C}_{i}+\mathbf{D}_{i}$ there exists the unique stationary probability vector $\widetilde{\theta }^{\left( i\right) }=\left(  \theta_{1}^{\left(  i\right)},\theta_{2}^{\left(  i\right)  },\cdots,\theta_{m}^{\left(  i\right)
}\right)  $  for $ 1\leq i\leq N$, that is, the vector $\widetilde{\theta }^{\left( i\right) }$ is the unique solution to the system of linear equations $\widetilde{\theta }^{\left( i\right) }\left( \mathbf{C}_{i}+\mathbf{D}_{i}\right) =0$ and $\widetilde{\theta }^{\left( i\right) }e=1$. In this case, the stationary average arrival rate of the MAP$\left( \mathbf{C}_{i}+\mathbf{D}_{i}\right) $ is $\lambda _{i}=\widetilde{\theta }^{\left( i\right) }D^{\left(
i\right) }e$. Specifically, we write that $\overrightarrow{\lambda }_{i}=\left( \lambda _{i}^{\left( 1\right)
},\lambda _{i}^{\left( 2\right) },\cdots ,\lambda _{i}^{\left( m\right)
}\right) =\widetilde{\theta }^{\left( i\right) }\mathbf{D}_{i}$ for $1\leq i\leq N$.

\textbf{(5) The first riding-bike time:} An outside bike user arrives at
the $i$th station to rent a bike. If there is no bike in the
$i$th station (i.e., the $i$th station is empty), then the user immediately
leaves this bike sharing system. If there is at least one available bike at the $i$th
station, then the user rents a bike and goes to Road $i\rightarrow
j$ for $j\neq i$ with probability
$p_{i,j}$ for $\sum_{j\in\Theta_{i}}p_{i,j}=1$, and his riding-bike time on Road $i\rightarrow j$
is an exponential random variable with riding-bike rate $\mu_{i,j}>0$.

\textbf{(6) The bike return times:}

Notice that for any user, his first bike return process may be different from those retrial processes with successively
returning the bike to one station for at least twice due to his pasted arrivals at the full
stations. In this situation, his road selection as well as his riding-bike
time in the first process may be different from those in any retrial return process.

\underline{The first return} -- When the user completes his short trip on
Road $i\rightarrow j$, he needs to return his bike to
the $j$th station. If there is at least one available parking position (i.e.,
a vacant docker), then the user directly returns the bike to the $j$th
station, and immediately leaves this bike sharing systems.

\underline{The second return} -- If no parking position is available at the
$j$th station, then the user has to ride the bike to the $l_{1}$th station with
probability $\alpha_{j,l_{1}}$ for $l_{1}\neq j$ and $\sum_{l_{1}\in\Theta
_{j}}\alpha_{j,l_{1}}=1$; and his future riding-bike time on Road $j\rightarrow l_{1}$ is also an
exponential random variable with riding-bike rate $\xi_{j,l_{1}}>0$. If there
is at least one available parking position, then the user directly returns
his bike to the $l_{1}$th station, and immediately leaves this bike sharing system.

\underline{The ($k+1$)st return for $k\geq2$ }-- We assume that this bike has
not been returned at any station yet through $k$ consecutive returns.
In this case, the user has to try his ($k+1$)st lucky return. Notice that
the user goes to the $l_{k}$th station from the $l_{k-1}$th full station
with probability $\alpha_{l_{k-1},l_{k}}$ for $l_{k}\neq l_{k-1}$ and
$\sum_{l_{k}\in\Theta_{l_{k-1}}}\alpha_{l_{k-1},l_{k}}=1$; and his riding-bike
time on Road $l_{k-1}\rightarrow l_{k}$ is an exponential random variable with
riding-bike rate $\xi_{l_{k-1},l_{k}}>0$. If there is at least one available
parking position, then the user directly returns his bike to the $l_{k}$th
station, and immediately leaves this bike sharing system; otherwise he has to
continuously ride his bike in order to try to return the bike to another station again.

We further assume that the returning-bike process is persistent in the sense
that the user must find a station with an empty position to return his
bike because the bike is a public property.

It is seen from the above description that the parameters: $p_{i,j}$ and
$\mu_{i,j}$, for $j\neq i$ and $1\leq i,j\leq N$, of the first return, may be
different from the parameters: $\alpha_{i,j}$ and $\xi_{i,j}$, for $j\neq i$
and $1\leq i,j\leq N$, of the $k$th return for $k\geq2$. This is due to a simple observation that the user
possibly deal with more things (for example, tourism, shopping, visiting friends and
so on) in the first return process, but he becomes only one
return task for returning his bike to one station during the $k$ successive return processes for $k\geq2$.

\textbf{(7) The departure discipline:} The user departure process has two
different cases: (a) An outside user directly leaves the bike sharing
system if he arrives at an empty station; and (b) if one user rents and
uses a bike, and he finally returns the bike to a station, then the user
completes his trip, and immediately leaves the bike sharing system.

We assume that all the above random variables are independent of each other. For
such a bike sharing system, Fig. 1 provides some intuitive physical interpretation for the bike sharing system.

\begin{figure}[ptb]
\setlength{\abovecaptionskip}{0.cm}  \setlength{\belowcaptionskip}{-0.cm}
\centering                           \includegraphics[width=9cm]{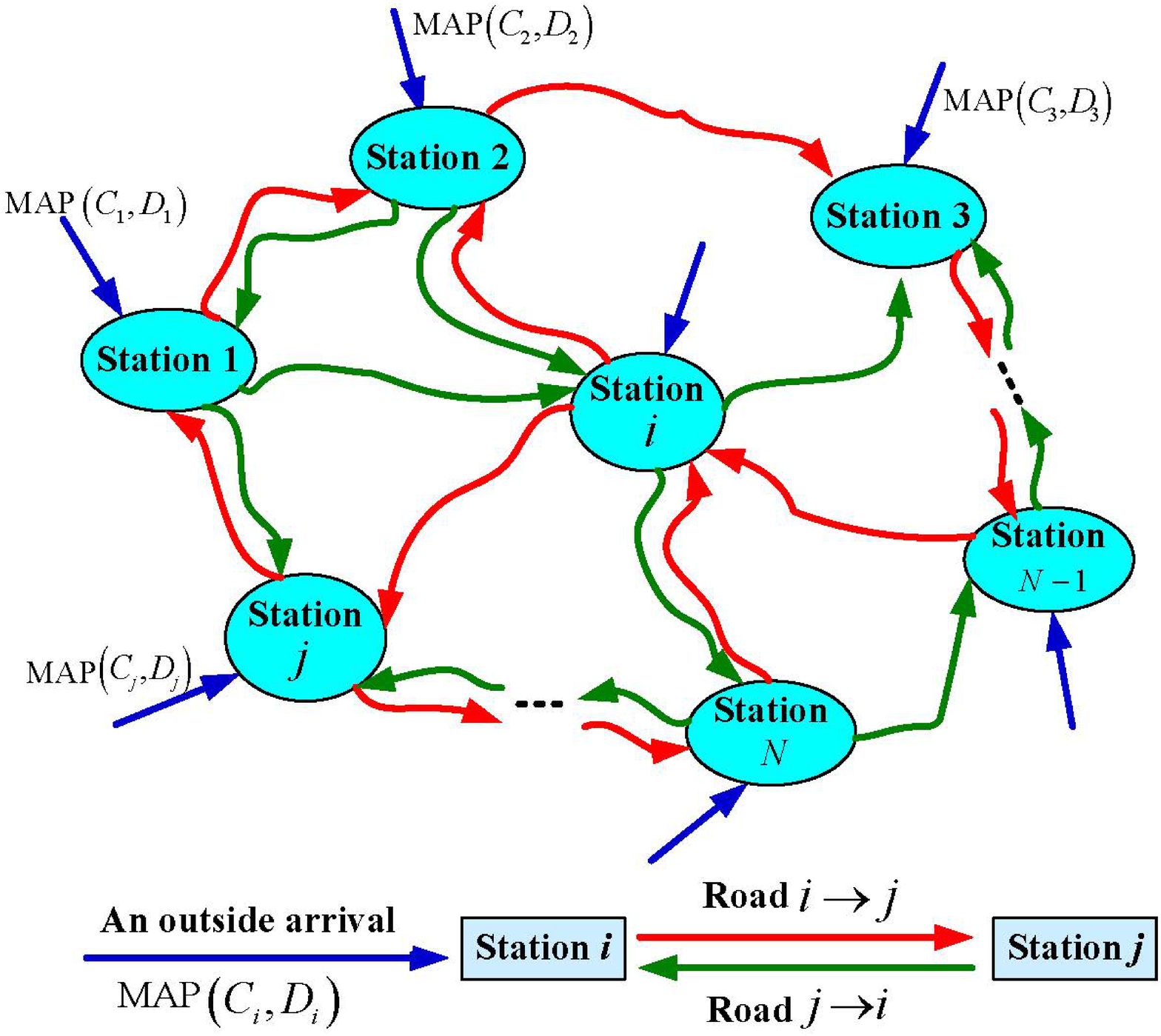}
\newline \caption{The physical structure of the bike sharing system}%
\label{figure:fig-1}%
\end{figure}

\section{A Closed Queueing Network}

In this section, we describe the bike sharing system as a closed queueing
network according to the fact that the number of bikes in this
system is fixed. To study such a closed queueing network, we need to determine the service
rates, the routing matrix and the relative arrival rates in all the virtual nodes.

For the bike sharing system, we need to abstract it as a closed queueing
network as follows:

\textbf{(1) Virtual nodes:} Although the stations and the roads have
different physical attributes, such as, different functions, different geographical topologies and so forth, it is seen that here the stations and the roads are all regarded as
the same abstructed nodes in a closed queueing network.

\textbf{(2) Virtual customers:} The bikes either at the stations or on the
roads are viewed as virtual customers as follows:

\underline{A closed queueing network under virtual idea:} The virtual customers are abstracted by the bikes from either the stations or the roads. In this case, the service processes are taken either from user arrivals at the station nodes or from users riding bikes on the road nodes. Since the
total number of bikes in the bike sharing system is fixed as the positive integer $NC$, thus the
bike sharing system can be regarded as a closed queueing network with such virtual customers, virtual nodes and virtual service times.

\underline{Two classes of virtual customers:} From Assumptions (2), (5) and (6) in
Section 2, it is seen that there are two different classes of virtual customers in the
road nodes, where the first class of virtual customers are the
bikes ridden on the roads for the first time; while the second class of
virtual customers are the bikes which are successively ridden on
the roads at least twice due to his arrivals at full stations.

We abstract the virtual nodes both from the stations and from the roads, and also find the virtual customers corresponding to the $NC$ bikes. This sets up a multiclass closed queueing network. To compute the steady-state probabilities of joint queue lengths in the bike
sharing system, it is seen from Chapter 7 in Bolch \emph{et al}. \cite{Bol:2006} that
we need to determine the service rate and the relative
arrival rate for each virtual node in the multiclass closed queueing network.

\textbf{(a) The service rates at nodes}

We discuss the service processes of the closed queueing network from two different cases: One for the station nodes, and the other for the road nodes. Fig. 2 shows how the two classes of service times are given from the multiclass closed queueing network.
\begin{figure}[ptb]
\setlength{\abovecaptionskip}{0.cm}  \setlength{\belowcaptionskip}{-0.cm}
\centering                            \includegraphics[width=10cm]{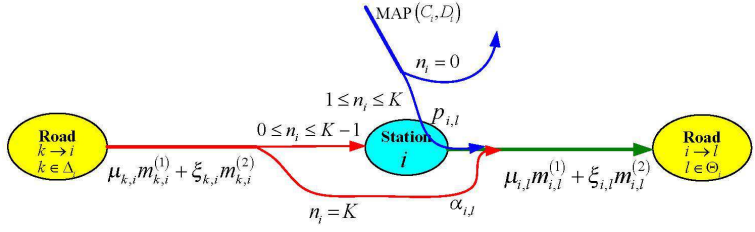}
\newline \caption{The queueing processes in the multiclass closed queueing network}%
\label{figure:fig-2}%
\end{figure}

\textbf{Case one:} A road node in the set $\cup_{i=1}^{N}R\left(
i\right)  $

\underline{The first class of virtual customers:} We denote the number of virtual customers of the first class on Road $i\rightarrow l$ by $m_{i,l}^{\left(  1\right)  }$. The return process of bikes of the first class from
Road $i\rightarrow l$ to Station $l$ for the first time is Poisson with service rate
\[¡¢
a_{i,l}^{\left(  1\right)  }=m_{i,l}^{\left(  1\right)  }\mu_{i,l}.
\]

\underline{The second class of virtual customers:}  We denote the number of virtual customers of the second class on Road $i\rightarrow l$ by $m_{i,l}^{\left(  2\right)  }$. The retrial return process of customers of the second class from Road $i\rightarrow j$ to Station $l$ is Poisson with service rate
\[
a_{i,l}^{\left(  2\right)  }=m_{i,l}^{\left(  2\right)  }\xi_{i,l}.
\]

\textbf{Case two:} The $N$ station nodes

Let $n_{i}$ be the number of bikes packed in Station $i$. The departure process of bikes from the $i$th station is due to those customers who rent the bikes at the $i$th station and then immediately enter one road in $R\left(  i\right)  $. Thus if the $i$th station is not empty, then the service process (i.e. renting bikes) is a MAP with a stationary service rate of phase $v$
\begin{equation}
a_{i}^{\left( v\right) }=\lambda _{i}^{\left( v\right) }\mathbf{1}_{\left\{
1\leq n_{i}\leq K\right\} }\sum_{l\neq i}^{N}p_{i,l}=\lambda _{i}^{\left(
v\right) }\mathbf{1}_{\left\{ 1\leq n_{i}\leq K\right\} },1\leq v\leq m,
\label{Cequ-5}
\end{equation}
where $\sum_{l\neq i}^{N}p_{i,l}=1$, and $\overrightarrow{\lambda }_{i}=\left( \lambda _{i}^{\left( 1\right)
},\lambda _{i}^{\left( 2\right) },\ldots ,\lambda _{i}^{\left( m\right)
}\right) $ is given by the MAP $ \left( C_{i},D_{i}\right) $ through $\overrightarrow{\lambda }_{i}=\widetilde{\theta }^{\left( i\right) }D_{i}$ for $1\leq i\leq N$.

\textbf{(b) The relative arrival rates}

For the multiclass closed queueing network, to determine the steady-state probability distribution of joint queue lengths at any  virtual node, it is necessary to firstly give the relative arrival rates at the virtual nodes. To this end, we must establish the routing matrix in the first step.

Based on Chapter 7 in Bolch \emph{et al}. \cite{Bol:2006}, we denote by $e_{i}$ and $e_{R_{i\rightarrow j}}^{\left( r\right) }$ the relative arrival rates of the $i$th station, and of Road $i\rightarrow l$ with bikes of class $r$, respectively. We
write%
\begin{equation*}
\overrightarrow{e}=\left\{ \overrightarrow{e}_{i}:1\leq i\leq N\right\} ,
\end{equation*}
where%
\begin{equation*}
\overrightarrow{e}_{i}=\left\{ \mathbf{e}_{i},\mathbf{e}_{R_{i\rightarrow
j}}^{\left( r\right) },j\in \Theta _{i},r=1,2\right\} .
\end{equation*}

Note that this bike sharing system is large-scale, thus the routing matrix of the closed queueing network corresponding to the bike sharing system will be very complicated. To understand how to set up such a routing matrix, in what follows we first give three simple examples for the purpose of writing the routing matrix, using the physical structure and the routing graph of the bike sharing system. See Figures 3 to 5 for more details.

Let $Q_{i}\left( t\right) $ be the number of bikes parked at Station $i$ at time $t\geq 0$. From the exponential and MAP assumptions, it is seen that an irreducible finite state Markov chain is used to express and analyze the bike sharing system, while the Markov chain is aperiodic and positive recurrent. In this case, there exists stationary probability vector in the Marokov chain, and thus we give the limit
\begin{equation*}
\pi _{i,K}=\underset{t\rightarrow +\infty }{\lim }P\left\{ Q_{i}\left(
t\right) =K\right\} .
\end{equation*}%

\textbf{Example One:} We consider a simple bike sharing system with two stations, and the physical structure of the stations and roads is depicted in (a) of Fig.3. Note that there exist two classes of virtual customers in the road nodes, and the bike routing graph of the bike sharing system is depicted in (b) of Fig.3. Since there are only two stations in this bike sharing system, we have $p_{i,j}=\alpha _{i,j}=1$. Based on this, we obtain the routing matrix of order 6 as follow:%
\begin{figure}[ptb]
\setlength{\abovecaptionskip}{0.cm}  \setlength{\belowcaptionskip}{-0.cm}
\centering
\includegraphics[width=12cm]{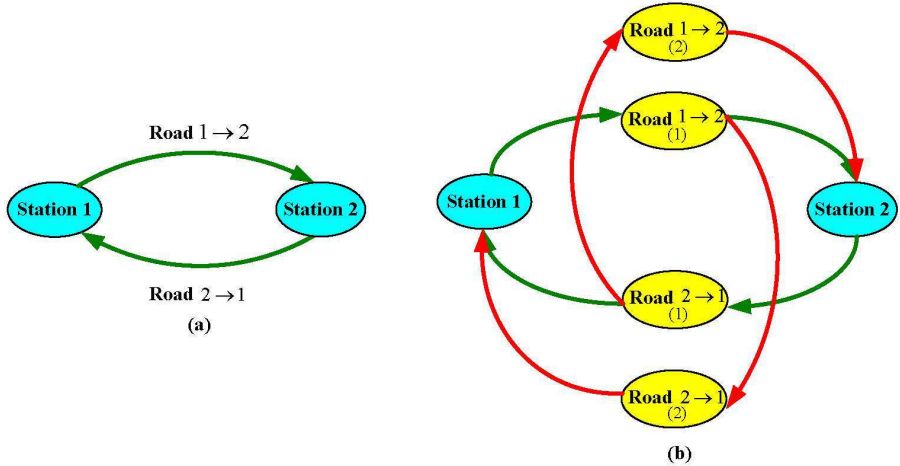}  \newline \caption{The physical
structure (a) and the bike routing graph (b) for a two-station bike sharing system}%
\label{figure:fig-3}%
\end{figure}

\[
P=\left(
\begin{array}
[c]{llllll}
& 1 &  &  &  & \\
&  &  & 1-\pi_{2,K} &  & \pi_{2,K}\\
&  &  & 1-\pi_{2,K} &  & \pi_{2,K}\\
&  &  &  & 1 & \\
1-\pi_{1,K} &  & \pi_{1,K} &  &  & \\
1-\pi_{1,K} &  & \pi_{1,K} &  &  &
\end{array}
\right)  ,
\]
where all those elements that are not expressed are viewed as zeros, and $\pi_{i,K}$ is a undetermined constant, and it is also the stationary probability of the $i$th full station for $i=1,2$.

To determine the relative arrival rate at each virtual node, using the system of linear equations $\overrightarrow{e}P=\overrightarrow{e}$ and $e_{1}=1$, we obtain%
\begin{equation*}
\left\{
\begin{array}{l}
e_{1}=\left( e_{R_{2\rightarrow 1}}^{\left( 1\right) }+e_{R_{2\rightarrow
1}}^{\left( 2\right) }\right) \left( 1-\pi _{1,K}\right) , \\
e_{R_{1\rightarrow 2}}^{\left( 1\right) }=e_{1}, \\
e_{R_{1\rightarrow 2}}^{\left( 2\right) }=\left( e_{R_{2\rightarrow
1}}^{\left( 1\right) }+e_{R_{2\rightarrow 1}}^{\left( 1\right) }\right) \pi
_{1,K}, \\
e_{2}=\left( e_{R_{1\rightarrow 2}}^{\left( 1\right) }+e_{R_{1\rightarrow
2}}^{\left( 2\right) }\right) \left( 1-\pi _{2,K}\right) , \\
e_{R_{2\rightarrow 1}}^{\left( 1\right) }=e_{2}, \\
e_{R_{2\rightarrow 1}}^{\left( 2\right) }=\left( e_{R_{1\rightarrow
2}}^{\left( 1\right) }+e_{R_{1\rightarrow 2}}^{\left( 2\right) }\right) \pi
_{2,K}.%
\end{array}%
\right.
\end{equation*}
Using $e_{1}=1$, we get%
\begin{equation}
\left\{
\begin{array}{l}
e_{1}=e_{R_{1\rightarrow 2}}^{\left( 1\right) }=1, \\
e_{R_{1\rightarrow 2}}^{\left( 2\right) }=\frac{\pi _{1,K}}{1-\pi _{1,K}},
\\
e_{2}=e_{R_{2\rightarrow 1}}^{\left( 1\right) }=\frac{1-\pi _{2,K}}{1-\pi
_{1,K}}, \\
e_{R_{2\rightarrow 1}}^{\left( 2\right) }=\frac{\pi _{2,K}}{1-\pi _{1,K}},%
\end{array}%
\right.   \label{sol-1}
\end{equation}%
where the two undetermined positive constants $\pi _{1,K}$ and $\pi _{2,K}$ will be given in the next section, and they determine the relative arrival rates at the six virtual nodes.

\textbf{Example Two:} We consider a bike sharing system with three stations, and the physical structure of the stations and roads can be seen in (a) of Fig.4. There exist two classes of virtual customers in the road nodes, and the bike routing graph of the bike sharing system is depicted in (b) of Fig.4. It is seen from (a) of Fig.4 that $ p_{1,2}=p_{2,3}=p_{3,1}=\alpha _{1,2}=\alpha _{2,3}=\alpha _{3,1}=1$. Based on this, the routing matrix of order 9 is given by%
\begin{figure}[ptb]
\setlength{\abovecaptionskip}{0.cm}  \setlength{\belowcaptionskip}{-0.cm}
\centering                           \includegraphics[width=6cm]{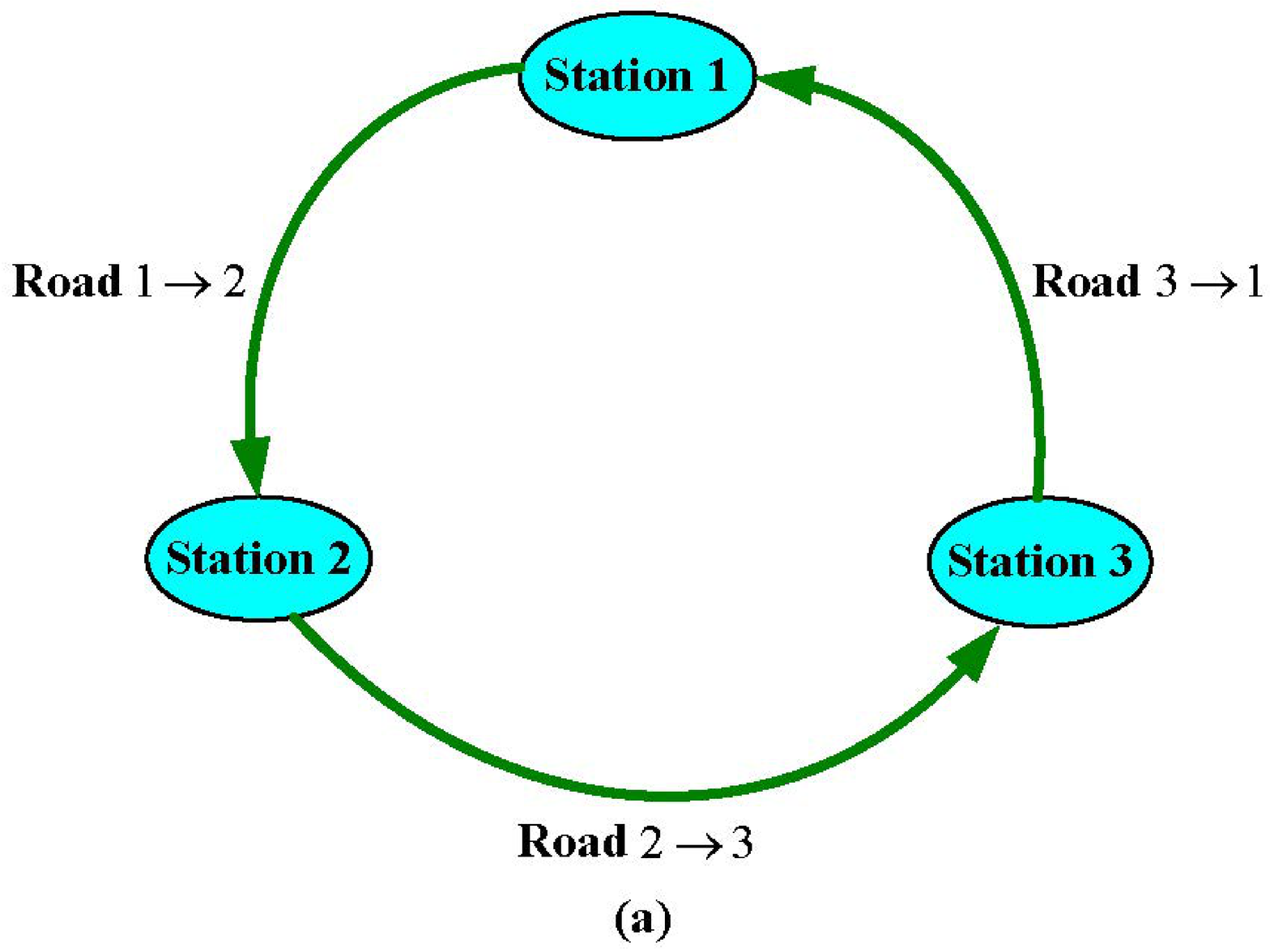}
\includegraphics[width=8cm]{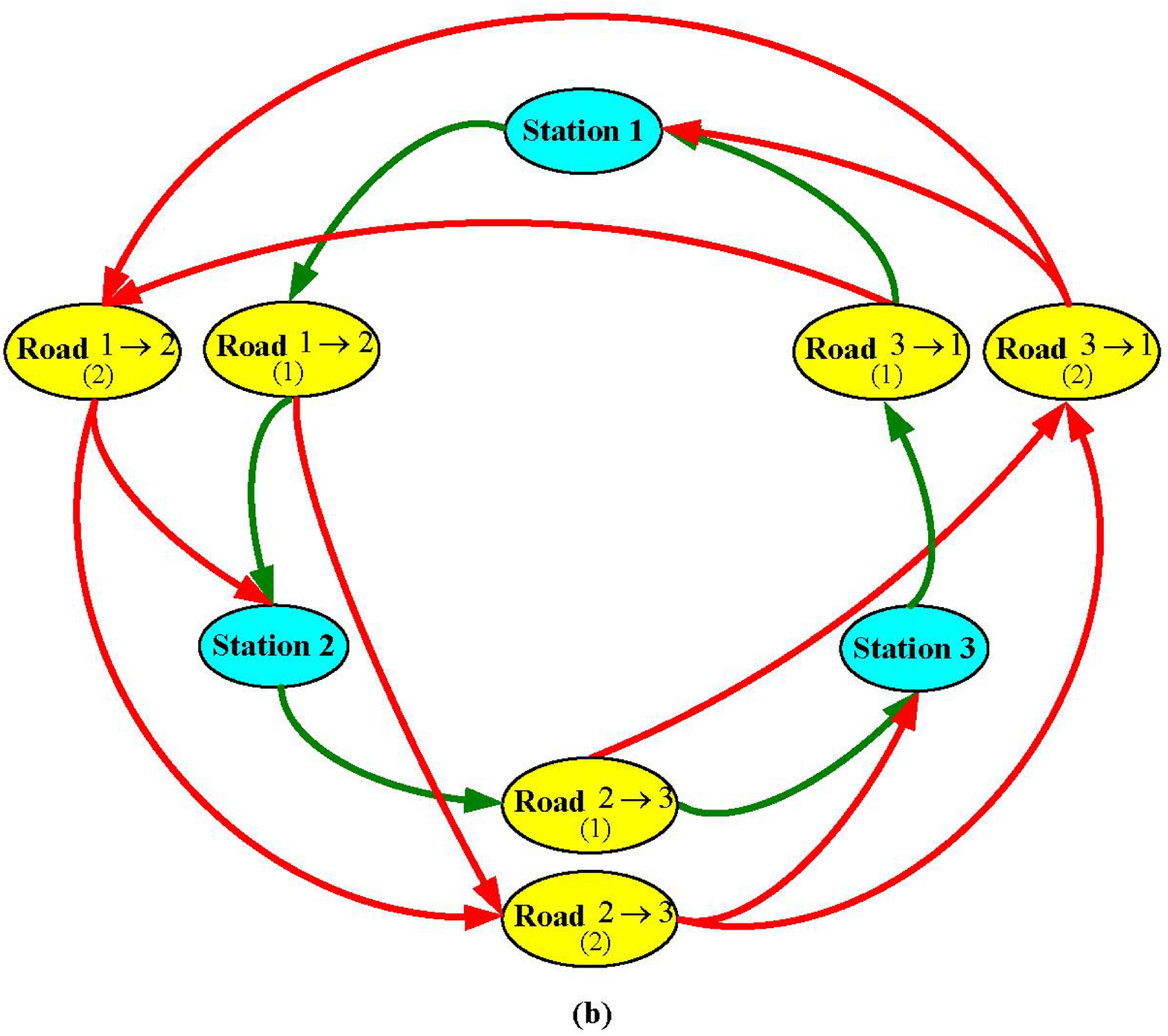}  \newline \caption{The physical
structure (a) and the bike routing graph (b) of a three-station bike sharing system}%
\label{figure:fig-4}%
\end{figure}

\[
\left(
\begin{array}
[c]{lllllllll}
& 1 &  &  &  &  &  &  & \\
&  &  & 1-\pi_{2,K} &  & \pi_{2,K} &  &  & \\
&  &  & 1-\pi_{2,K} &  & \pi_{2,K} &  &  & \\
&  &  &  & 1 &  &  &  & \\
&  &  &  &  &  & 1-\pi_{3,K} &  & \pi_{3,K}\\
&  &  &  &  &  & 1-\pi_{3,K} &  & \pi_{3,K}\\
&  &  &  &  &  &  & 1 & \\
1-\pi_{1,K} &  & \pi_{1,K} &  &  &  &  &  & \\
1-\pi_{1,K} &  & \pi_{1,K} &  &  &  &  &  &
\end{array}
\right)  .
\]

To determine the relative arrival rate at each virtual node, using the system of linear equations $\overrightarrow{e}P=\overrightarrow{e}$ and $e_{1}=1$, we obtain%
\begin{equation*}
\left\{
\begin{array}{l}
e_{1}=\left( e_{R_{3\rightarrow 1}}^{\left( 1\right) }+e_{R_{3\rightarrow
1}}^{\left( 2\right) }\right) \left( 1-\pi _{1,K}\right) , \\
e_{R_{1\rightarrow 2}}^{\left( 1\right) }=e_{1}, \\
e_{R_{1\rightarrow 2}}^{\left( 2\right) }=\left( e_{R_{3\rightarrow
1}}^{\left( 1\right) }+e_{R_{3\rightarrow 1}}^{\left( 2\right) }\right) \pi
_{1,K}, \\
e_{2}=\left( e_{R_{1\rightarrow 2}}^{\left( 1\right) }+e_{R_{1\rightarrow
2}}^{\left( 2\right) }\right) \left( 1-\pi _{2,K}\right) , \\
e_{R_{2\rightarrow 3}}^{\left( 1\right) }=e_{2}, \\
e_{R_{2\rightarrow 3}}^{\left( 2\right) }=\left( e_{R_{1\rightarrow
2}}^{\left( 1\right) }+e_{R_{1\rightarrow 2}}^{\left( 2\right) }\right) \pi
_{2,K}, \\
e_{3}=\left( e_{R_{2\rightarrow 3}}^{\left( 1\right) }+e_{R_{2\rightarrow
3}}^{\left( 2\right) }\right) \left( 1-\pi _{3,K}\right) , \\
e_{R_{3\rightarrow 1}}^{\left( 1\right) }=e_{3}, \\
e_{R_{3\rightarrow 1}}^{\left( 2\right) }=\left( e_{R_{2\rightarrow
3}}^{\left( 1\right) }+e_{R_{2\rightarrow 3}}^{\left( 2\right) }\right) \pi
_{3,K}.%
\end{array}%
\right.
\end{equation*}%
Using $e_{1}=1$, we get%
\begin{equation}
\left\{
\begin{array}{l}
e_{1}=e_{R_{1\rightarrow 2}}^{\left( 1\right) }=1, \\
e_{R_{1\rightarrow 2}}^{\left( 2\right) }=\frac{\pi _{1,K}}{1-\pi _{1,K}} \\
e_{2}=e_{R_{2\rightarrow 3}}^{\left( 1\right) }=\frac{1-\pi _{2,K}}{1-\pi
_{1,K}}, \\
e_{R_{2\rightarrow 3}}^{\left( 2\right) }=\frac{\pi _{2,K}}{1-\pi _{1,K}},
\\
e_{3}=e_{R_{3\rightarrow 1}}^{\left( 1\right) }=\frac{1-\pi _{3,K}}{1-\pi
_{1,K}}, \\
e_{R_{3\rightarrow 1}}^{\left( 2\right) }=\frac{\pi _{3,K}}{1-\pi _{1,K}},%
\end{array}%
\right.   \label{sol-2}
\end{equation}%
where the three undetermined positive constants $\pi _{1,K}$, $\pi _{2,K}$ and $\pi _{3,K}$ will be given in the next section, and they determine the relative arrival rates for the nine virtual nodes.

\textbf{Example Three:} We consider a bike sharing system with three stations, and the physical structure of the stations and roads can be seen in (a) of Fig.5. There exist two classes of virtual customers in the road nodes, and the bike routing graph of the bike sharing system is depicted in (b) of Fig.5. Based on this, we obtain the routing matrix of order 11 as follow:%
\begin{figure}[ptb]
\setlength{\abovecaptionskip}{0.cm}  \setlength{\belowcaptionskip}{-0.cm}
\centering                           \includegraphics[width=5cm]{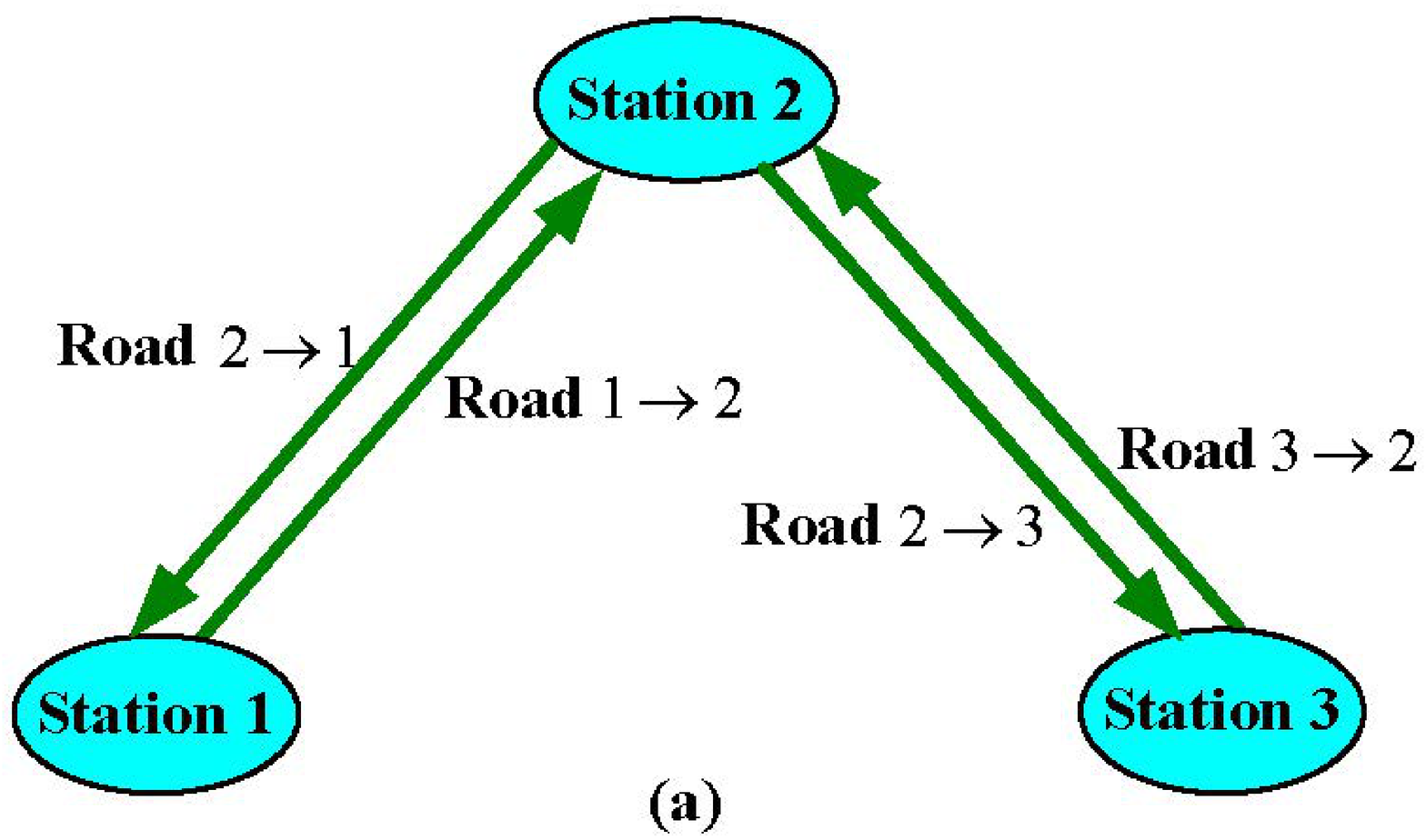}
\includegraphics[width=10cm]{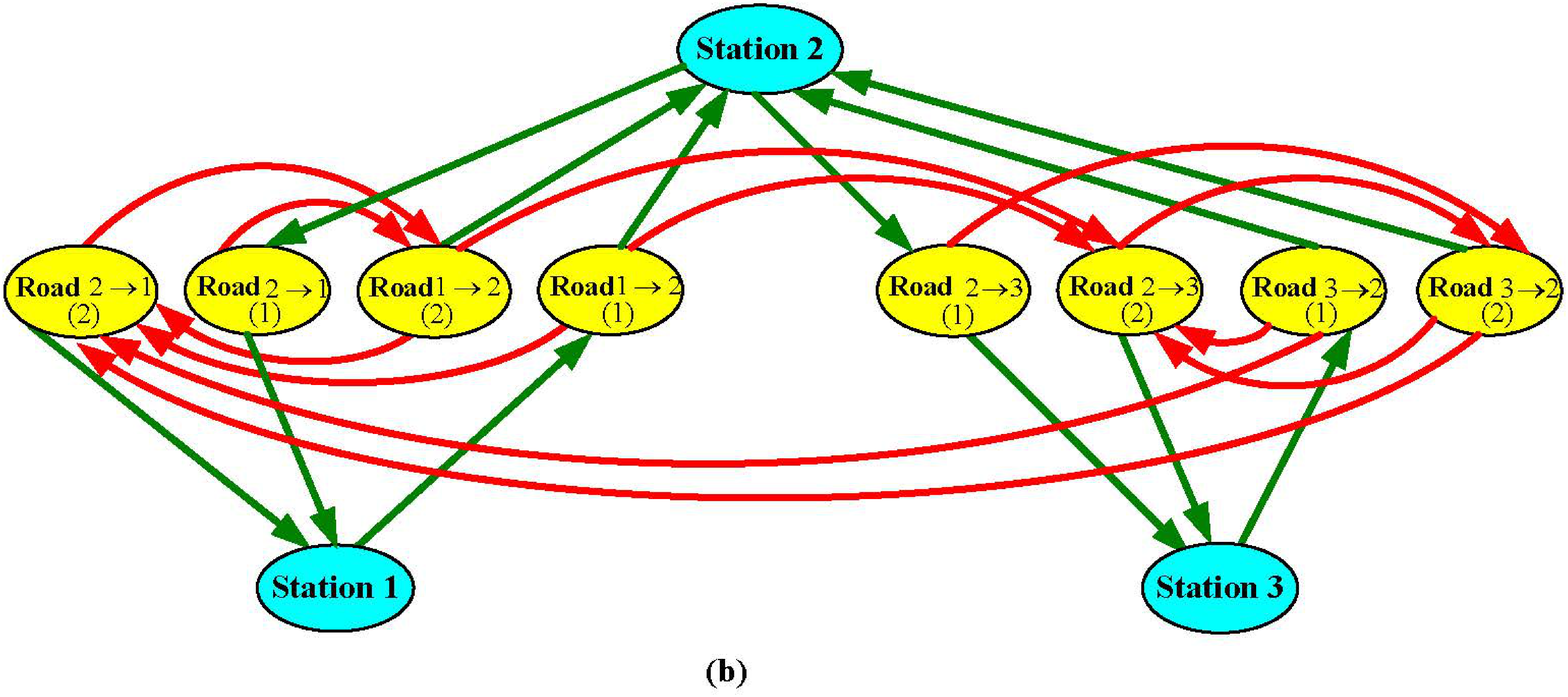}  \newline \caption{The physical
structure (a) and the bike routing graph (b) of a three-station bike sharing system}%
\label{figure:fig-5}%
\end{figure}
\begin{equation*}
\left(
\begin{array}{ccccccccccc}
& 1 &  &  &  &  &  &  &  &  &  \\
&  &  & 1-\pi _{2,K} &  & \alpha _{2,1}\pi _{2,K} &  & \alpha _{2,3}\pi
_{2,K} &  &  &  \\
&  &  & 1-\pi _{2,K} &  & \alpha _{2,1}\pi _{2,K} &  & \alpha _{2,3}\pi
_{2,K} &  &  &  \\
&  &  &  & p_{2,1} &  & p_{2,3} &  &  &  &  \\
1-\pi _{1,K} &  & \pi _{1,K} &  &  &  &  &  &  &  &  \\
1-\pi _{1,K} &  & \pi _{1,K} &  &  &  &  &  &  &  &  \\
&  &  &  &  &  &  &  & 1-\pi _{3,K} &  & \pi _{3,K} \\
&  &  &  &  &  &  &  & 1-\pi _{3,K} &  & \pi _{3,K} \\
&  &  &  &  &  &  &  &  & 1 &  \\
&  &  & 1-\pi _{2,K} &  & \alpha _{2,1}\pi _{2,K} &  & \alpha _{2,3}\pi
_{2,K} &  &  &  \\
&  &  & 1-\pi _{2,K} &  & \alpha _{2,1}\pi _{2,K} &  & \alpha _{2,3}\pi
_{2,K} &  &  &
\end{array}%
\right) .
\end{equation*}
To determine the relative arrival rate at each virtual node, using the system of linear equations $\overrightarrow{e}P=\overrightarrow{e}$ and $e_{1}=1$, we obtain%
\begin{equation}
\left\{
\begin{array}{l}
e_{1}=\left( e_{R_{2\rightarrow 1}}^{\left( 1\right) }+e_{R_{2\rightarrow
1}}^{\left( 2\right) }\right) \left( 1-\pi _{1,K}\right) , \\
e_{R_{1\rightarrow 2}}^{\left( 1\right) }=e_{1}, \\
e_{R_{1\rightarrow 2}}^{\left( 2\right) }=\left( e_{R_{2\rightarrow
1}}^{\left( 1\right) }+e_{R_{2\rightarrow 1}}^{\left( 2\right) }\right) \pi
_{1,K} \\
e_{2}=\left( e_{R_{1\rightarrow 2}}^{\left( 1\right) }+e_{R_{1\rightarrow
2}}^{\left( 2\right) }+e_{R_{3\rightarrow 2}}^{\left( 1\right)
}+e_{R_{3\rightarrow 2}}^{\left( 2\right) }\right) \left( 1-\pi
_{2,K}\right) , \\
e_{R_{2\rightarrow 1}}^{\left( 1\right) }=p_{2,1}e_{2}, \\
e_{R_{2\rightarrow 1}}^{\left( 2\right) }=\left( e_{R_{1\rightarrow
2}}^{\left( 1\right) }+e_{R_{1\rightarrow 2}}^{\left( 2\right)
}+e_{R_{3\rightarrow 2}}^{\left( 1\right) }+e_{R_{3\rightarrow 2}}^{\left(
2\right) }\right) \alpha _{2,1}\pi _{2,K}, \\
e_{R_{2\rightarrow 3}}^{\left( 1\right) }=p_{2,3}e_{2}, \\
e_{R_{2\rightarrow 3}}^{\left( 2\right) }=\left( e_{R_{1\rightarrow
2}}^{\left( 1\right) }+e_{R_{1\rightarrow 2}}^{\left( 2\right)
}+e_{R_{3\rightarrow 2}}^{\left( 1\right) }+e_{R_{3\rightarrow 2}}^{\left(
2\right) }\right) \alpha _{2,3}\pi _{2,K}, \\
e_{3}=\left( e_{R_{2\rightarrow 3}}^{\left( 1\right) }+e_{R_{2\rightarrow
3}}^{\left( 2\right) }\right) \left( 1-\pi _{3,K}\right) , \\
e_{R_{3\rightarrow 2}}^{\left( 1\right) }=e_{3}, \\
e_{R_{3\rightarrow 2}}^{\left( 2\right) }=\left( e_{R_{2\rightarrow
3}}^{\left( 1\right) }+e_{R_{2\rightarrow 3}}^{\left( 2\right) }\right) \pi
_{3,K}.%
\end{array}%
\right.   \label{e-3}
\end{equation}%
By using $e_{1}=1$, we obtain%
\begin{equation*}
\left\{
\begin{array}{l}
e_{1}=e_{R_{1\rightarrow 2}}^{\left( 1\right) }=1, \\
e_{R_{1\rightarrow 2}}^{\left( 2\right) }=\frac{\pi _{1,K}}{1-\pi _{1,K}} \\
e_{2}=\frac{1-\pi _{2,K}}{\left( 1-\pi _{1,K}\right) \left[ \left( \alpha
_{2,1}-p_{2,1}\right) \pi _{2,K}+p_{2,1}\right] }, \\
e_{R_{2\rightarrow 1}}^{\left( 1\right) }=\frac{p_{2,1}\left( 1-\pi
_{2,K}\right) }{\left( 1-\pi _{1,K}\right) \left[ \left( \alpha
_{2,1}-p_{2,1}\right) \pi _{2,K}+p_{2,1}\right] }, \\
e_{R_{2\rightarrow 1}}^{\left( 2\right) }=\frac{\alpha _{2,1}\pi _{2,K}}{%
\left( 1-\pi _{1,K}\right) \left[ \left( \alpha _{2,1}-p_{2,1}\right) \pi
_{2,K}+p_{2,1}\right] }, \\
e_{R_{2\rightarrow 3}}^{\left( 1\right) }=\frac{p_{2,3}\left( 1-\pi
_{2,K}\right) }{\left( 1-\pi _{1,K}\right) \left[ \left( \alpha
_{2,1}-p_{2,1}\right) \pi _{2,K}+p_{2,1}\right] }, \\
e_{R_{2\rightarrow 3}}^{\left( 2\right) }=\frac{\alpha _{2,3}\pi _{2,K}}{%
\left( 1-\pi _{1,K}\right) \left[ \left( \alpha _{2,1}-p_{2,1}\right) \pi
_{2,K}+p_{2,1}\right] }, \\
e_{3}=e_{R_{3\rightarrow 1}}^{\left( 1\right) }=\frac{\left( 1-\pi
_{3,K}\right) \left[ \left( \alpha _{2,3}-p_{2,3}\right) \pi _{2,K}+p_{2,3}%
\right] }{\left( 1-\pi _{1,K}\right) \left[ \left( \alpha
_{2,1}-p_{2,1}\right) \pi _{2,K}+p_{2,1}\right] }, \\
e_{R_{3\rightarrow 1}}^{\left( 2\right) }=\frac{\pi _{3,K}\left[ \left(
\alpha _{2,3}-p_{2,3}\right) \pi _{2,K}+p_{2,3}\right] }{\left( 1-\pi
_{1,K}\right) \left[ \left( \alpha _{2,1}-p_{2,1}\right) \pi _{2,K}+p_{2,1}%
\right] }.%
\end{array}%
\right.
\end{equation*}%
\\ [0.5cm]
\textbf{The routing matrices for more general case}

Observing the three examples, it may be easy and convenient to write a routing matrix for a more general bike sharing system. Note that Example Three provides more intuitive understanding on how to write those elements of the routing matrix, thus for a more general bike sharing system we establish the routing matrix $ \mathbf{P}=\left( p_{\widetilde{i},\widetilde{j}}\right) $ as follow:
\begin{equation*}
p_{\widetilde{i},\widetilde{j}}=\left\{
\begin{array}{ll}
p_{i,j}, & \text{if }\widetilde{i}=\text{Station }i, \widetilde{j}=\text{Road
}i\rightarrow j \\
1-\pi _{j,K}, & \text{if }\widetilde{i}=\text{Road }i\rightarrow j,
\widetilde{j}=\text{Station }j \\
\sum\limits_{l\in \Theta _{i}\&l\in \Delta _{j}}\alpha _{l,j}\pi _{l,K}, &
\text{if }n_{l}=K, \widetilde{i}=\text{Road }i\rightarrow l, \widetilde{j}=%
\text{Road }l\rightarrow j \\
0. & \text{otherwise}%
\end{array}%
\right.
\end{equation*}
\begin{The}
The routing matrix $\mathbf{P}$ of finite size is irreducible and stochastic, and
there exists the unique positive solution to the following system of linear
equations%
\begin{equation*}
\left\{
\begin{array}{l}
\overrightarrow{e}=\overrightarrow{e}\mathbf{P}, \\
e_{1}=1,%
\end{array}%
\right.
\end{equation*}%
where $e_{1}=1$ is the first element of the
row vector $\overrightarrow{e}$, and $\overrightarrow{e}$ is a row vector of the relative arrival rates of this bike sharing system.
\end{The}

\textbf{Proof:} The outline of this proof is described as follows. It is clear that the size of the routing matrix $\mathbf{P}$ is finite. At the same time, it is well-known that (a) the routing structure of the multiclass closed queueing network indicates
that the routing matrix $\mathbf{P}$ is stochastic; and (b) the accessibility of
each station node or road node in the bike sharing system shows that the routing matrix
$\mathbf{P}$ is irreducible. Thus the routing matrix $\mathbf{P}$
is not only irreducible but also stochastic. For the routing matrix
$\mathbf{P}$, applying Theorem 1.1 (a) and (b) of Chapter 1 in Seneta
\cite{Sen:1981}, the left eigenvector $\overrightarrow{e}$ of the
irreducible stochastic matrix $\mathbf{P}$ of finite sizes corresponding to the maximal
eigenvalue $1$ is strictly positive, that is, $\overrightarrow{e}>0$; and
$\overrightarrow{e}$ is unique with $e_{1}=1$. This completes this proof. \textbf{{\rule{0.08in}{0.08in}}}

\textbf{(c) A joint queue-length process}

Let $Q_{i}^{\left(  v\right)  }\left(
t\right)  $\ be the number of bikes parked in Station $i$ with phase $v$ of the
MAP at time $t\geq0$, for $1\leq i\leq N$, $1\leq v\leq m$; and $R_{k,l}^{\left( r\right) }\left( t\right) $ the number of bikes of class $r$
ridden on Road $k\rightarrow l$ at time $t\geq0$, for $r=1,2$ and for $l\neq k$
with $1\leq k,l\leq N$. We write%
\[
\mathbf{X}\left(  t\right)  =\left(  \mathbf{L}_{1}\left(  t\right)
,\mathbf{L}_{2}\left(  t\right)  ,\ldots,\mathbf{L}_{N-1}\left(  t\right)
,\mathbf{L}_{N}\left(  t\right)  \right)  ,
\]
where for $1\leq i\leq N$%
\begin{equation*}
\mathbf{L}_{i}\left( t\right) =\left( Q_{i}^{\left( 1\right) }\left(
t\right) ,Q_{i}^{\left( 2\right) }\left( t\right) ,\ldots ,Q_{i}^{\left(
m\right) }\left( t\right) ;R_{i,j}^{\left( 1\right) }\left( t\right)
,R_{i,j}^{\left( 2\right) }\left( t\right) ,\text{ }j\in \Theta _{i}\right) .
\end{equation*}
Obviously, $\left\{  \mathbf{X}\left(  t\right)  :t\geq0\right\}  $ is a
Markov process due to the exponential and MAP assumptions of this bike sharing system. It is easy to see that the state space of Markov
process $\left\{  \mathbf{X}\left(  t\right)  :t\geq0\right\}  $ is given by%
\begin{align}
\Omega  =&\left\{ \overrightarrow{n}:0\leq n_{i}^{\left( v\right) }\leq K,%
\text{ }1\leq i\leq N,1\leq v\leq m,\right.   \notag \\
& \left. 0\leq m_{k,l}^{\left( r\right) }\leq NC,\text{ }r=1,2,l\neq k,1\leq
k,l\leq N,\right.   \label{o-1} \\
& \left. \sum\limits_{i=1}^{N}\sum\limits_{v=1}^{m}n_{i}^{\left( v\right)
}\left( t\right) +\sum_{k=1}^{N}\sum_{l\in \Theta
_{k}}\sum_{r=1}^{2}m_{k,l}^{\left( r\right) }=NC\right\} ,  \notag
\end{align}
where%
\begin{equation*}
\overrightarrow{n}=\left( \mathbf{n}_{1},\mathbf{n}_{2},\ldots ,\mathbf{n}%
_{N}\right) ,
\end{equation*}%
for $1\leq i\leq N$%
\begin{equation*}
\mathbf{n}_{i}=\left( n_{i}^{\left( 1\right) },n_{i}^{\left( 2\right)
},\cdots ,n_{i}^{\left( m\right) };m_{i,j}^{\left(
1\right) },m_{i,j}^{\left( 2\right) },j\in \Theta _{i}\right) .
\end{equation*}%

It is easy to check that the Markov process $\left\{ \mathbf{X}\left( t\right) :t\geq 0\right\} $ on a finite state space is irreducible, aperiodic and positive recurrent. Therefore, there exists the stationary probability vector
\begin{equation*}
\pi =\left( \pi \left( \overrightarrow{n}\right) :\overrightarrow{n}\in
\Omega \right)
\end{equation*}%
such that%
\begin{equation*}
\pi \left( \overrightarrow{n}\right) =\underset{t\rightarrow +\infty }{\lim }%
P\left\{ \mathbf{X}\left( t\right) =\overrightarrow{n}\right\} .
\end{equation*}

\section{A Product-Form Solution and Performance Analysis}

In this section, we first provide a product-form solution to the steady-state
probabilities of joint queue lengths in the multiclass closed queueing network. Then we provide
a nonlinear solution to determine the $N$ undetermined constants: $\pi _{1,K},\pi _{2,K},\ldots ,\pi _{N,K}$. Also, an example is used to indicate our computational steps. Finally, we analyze performance measures of the bike sharing system by means of the steady-state
probabilities of joint queue lengths.

Note that $\left\{  \mathbf{X}\left(  t\right)  :t\geq0\right\}  $ is an
irreducible, aperiodic, positive recurrent and continuous-time Markov process with finite states, thus we have
\begin{align*}
\mathbf{\pi }\left( \overrightarrow{n}\right) =& \lim_{t\rightarrow +\infty
}P\left\{ Q_{i}^{\left( v\right) }\left( t\right) =n_{i}^{\left( v\right)
},1\leq i\leq N,1\leq v\leq m;\text{ }R_{k,l}^{\left( 1\right) }\left(
t\right) =m_{k,l}^{\left( 1\right) },R_{k,l}^{\left( 2\right) }\left(
t\right) \right.  \\
& \left. =m_{k,l}^{\left( 2\right) },1\leq k,l\leq N\text{ with }l\neq
k,\sum\limits_{i=1}^{N}\sum\limits_{v=1}^{m}n_{i}^{\left( v\right)
}+\sum\limits_{k=1}^{N}\sum\limits_{l\in \Theta
_{k}}\sum_{r=1,2}m_{k,l}^{\left( r\right) }=NC\right\} .
\end{align*}
Note that if $\sum\nolimits_{i=1}^{N}\sum\nolimits_{v=1}^{m}n_{i}^{\left(
v\right) }+\sum\nolimits_{k=1}^{N}\sum\nolimits_{l\in \Theta
_{k}}\sum_{r=1,2}m_{k,l}^{\left( r\right) }\neq NC$, it is easy to see that $\pi \left( \overrightarrow{n}\right) =0$.
In practice, it is a key in the study of bike sharing systems to provide expression for the steady-state probability $\pi \left( \overrightarrow{n}\right) $, $\overrightarrow{n}\in \Omega $.

\subsection{A product-form solution}

For the bike sharing system, we establish a multiclass closed queueing network with $N+\sum\nolimits_{i=1}^{N}\left| R\left( i\right) \right| $ virtual nodes and with $NC$ virtual customers. As $t\rightarrow +\infty $, the multiclass closed queueing network is decomposed into $N+\sum\nolimits_{i=1}^{N}\left| R\left( i\right) \right| $ isolated and equivalent queueing systems as follows:

(i) \textbf{The $i$th station node:} An equivalent queue is M$_{i}$/MAP$_{i}$/1/K, where M$_{i}$ denotes a Poisson process with relative arrival rate $e_{i}$, and MAP$_{i}$ is MAP$\left( C_{i},D_{i}\right) $ as a service process.

(ii) \textbf{The Road $i\rightarrow l$ node:} The two classes of customers correspond to their two queueing processes as follow:

(a) The first queue process on the Road $i\rightarrow l$ node is M$_{i\rightarrow j}^{\left( 1\right) }$/$%
\sum\nolimits_{k=1}^{m_{i,j}^{\left( 1\right) }}$M$_{i\rightarrow
j;1}^{\left( k\right) }$/1, where M$_{i\rightarrow j}^{\left( 1\right) }$ denotes
a Poisson process with relative arrival rate $e_{R_{i\rightarrow j}}^{\left(
1\right) }$, and $\sum\nolimits_{k=1}^{m_{i,j}^{\left( 1\right) }}$M$%
_{i\rightarrow j;1}^{\left( k\right) }$ is the random sum of $m_{i,j}^{\left( 1\right) }$ i.i.d. exponential random variables, each of which is exponential with service rate $%
\mu _{i,j}$.

(b) The second queue process on the Road $i\rightarrow l$ node is M$_{i\rightarrow j}^{\left( 2\right) }$/$%
\sum\nolimits_{k=1}^{m_{i,j}^{\left( 2\right) }}$M$_{i\rightarrow
j;2}^{\left( k\right) }$/1, in which M$_{i\rightarrow j}^{\left( 2\right) }$
is a Poisson process with relative arrival rate $e_{R_{i\rightarrow
j}}^{\left( 2\right) }$, and $\sum\nolimits_{k=1}^{m_{i,j}^{\left( 2\right) }}$M$%
_{i\rightarrow j;2}^{\left( k\right) }$ is the random sum of $m_{i,j}^{\left( 2\right) }$ i.i.d. exponential random variables, each of which is exponential with service rate $\xi _{i,j}$.

Using the above three classes of isolated queues, the following theorem provides a product-form solution to the steady-state probability $\mathbf{\pi}\left(  \overrightarrow{n}\right)  $ of joint queue lengths at the virtual nodes for
$\overrightarrow{n}\in\Omega$; while its proof is easy by means of Chapter 7
in Bolch \emph{et al}. \cite{Bol:2006} and is omitted here.

\begin{The}
For the two-class closed queueing network corresponding to the bike sharing system, if the undetermined constants $\pi _{1,K},\pi _{2,K},\ldots ,\pi _{N,K}$ are given, then
the steady-state
joint probability $\mathbf{\pi}\left(  \overrightarrow{n}\right)  $ is given
by
\begin{equation}
\mathbf{\pi }\left( \overrightarrow{n}\right) =\frac{1}{G\left( NC\right) }%
\prod_{i=1}^{N}H\left( \mathbf{n}_{i}\right) H\left( \mathbf{m}_{i}\right) ,
\label{Pro-1}
\end{equation}%
where $\overrightarrow{n}\in\Omega$,
\begin{equation*}
H\left( \mathbf{n}_{i}\right) =\frac{\left( n_{i}^{\left( 1\right)
}+n_{i}^{\left( 2\right) }+\cdots +n_{i}^{\left( m\right) }\right) !}{%
n_{i}^{\left( 1\right) }!n_{i}^{\left( 2\right) }!\cdots n_{i}^{\left(
m\right) }!}\prod_{v=1}^{m}\left( \frac{e_{i}}{\lambda _{i}^{\left( v\right)
}}\right) ^{n_{i}^{\left( v\right) }},
\end{equation*}%
\begin{equation*}
H\left( \mathbf{m}_{i}\right) =\prod_{j\in \Theta _{i}}\frac{\left(
m_{i,j}^{\left( 1\right) }+m_{i,j}^{\left( 2\right) }\right) !}{%
m_{i,j}^{\left( 1\right) }!m_{i,j}^{\left( 2\right) }!}\left( \frac{%
e_{R_{i\rightarrow j}}^{\left( 1\right) }}{m_{i,j}^{\left( 1\right) }\mu
_{i,j}}\right) ^{m_{i,j}^{\left( 1\right) }}\left( \frac{e_{R_{i\rightarrow
j}}^{\left( 2\right) }}{m_{i,j}^{\left( 2\right) }\xi _{i,j}}\right)
^{m_{i,j}^{\left( 2\right) }},
\end{equation*}%
and $G\left(  NC\right)  $ is a normalization constant, given by%
\begin{equation*}
G\left( NC\right) =\sum\limits_{\overrightarrow{n}\in \Omega
}\prod_{i=1}^{N}H\left( \mathbf{n}_{i}\right) H\left( \mathbf{m}_{i}\right) .
\end{equation*}%
\end{The}

By means of the product-form solution given in Theorem 2, the following theorem further establishes a system of nonlinear equations, whose solution determines the $N$ undetermined constants $\pi _{1,K},\pi _{2,K},\ldots ,\pi _{N,K}$. Note that $\pi _{i,K}$ is also the steady-state probability of the $i$th full station for $1\leq i\leq N$. While its proof is easy by means of the law of total probability and is omitted here.

\begin{The}
The undetermined constants $\pi _{1,K}$, $\pi _{2,K}$,$\ldots $, $\pi _{N,K}$ can be uniquely determined by the following system of nonlinear equations:
\begin{equation*}
\left\{
\begin{array}{c}
\pi _{1,K}=\sum\limits_{\substack{ \overrightarrow{n}\in \Omega  \\ %
\&n_{1}=K,}}\mathbf{\pi }\left( \overrightarrow{n}\right) , \\
\pi _{2,K}=\sum\limits_{\substack{ \overrightarrow{n}\in \Omega  \\ %
\&n_{2}=K,}}\mathbf{\pi }\left( \overrightarrow{n}\right) , \\
\vdots  \\
\pi _{N,K}=\sum\limits_{\substack{ \overrightarrow{n}\in \Omega  \\ %
\&n_{N}=K,}}\mathbf{\pi }\left( \overrightarrow{n}\right) ,%
\end{array}%
\right.
\end{equation*}
where $\pi \left( \overrightarrow{n}\right) $ is given by the product-form solution stated in Theorem 2.
\end{The}

To indicate how to compute the undetermined constants $\pi _{1,K},\pi _{2,K},\ldots ,\pi _{N,K}$, in what follows we give a concrete example.

\textbf{Example Four:} In Example One, we use the product-form solution to determine $\pi _{1,K}$ and $\pi _{2,K}$. By using (\ref{sol-1}) and (\ref{Pro-1}), we
obtain%
\begin{equation}
\left\{
\begin{array}{c}
\pi _{1,K}=\sum\limits_{\substack{ \overrightarrow{n}\in \Omega  \\ %
\&n_{1}=K,}}\mathbf{\pi }\left( \overrightarrow{n}\right) , \\
\pi _{2,K}=\sum\limits_{\substack{ \overrightarrow{n}\in \Omega  \\ %
\&n_{2}=K,}}\mathbf{\pi }\left( \overrightarrow{n}\right) .%
\end{array}%
\right.   \label{eqp-1}
\end{equation}%
We take that $C=2, K=3, m=2$. Thus (\ref{eqp-1}) is simplified as
\begin{equation}
\left\{
\begin{array}{l}
\pi _{1,K}=\frac{1}{G(NC)}\left( \frac{1}{\lambda _{1}^{\left( 1\right) }}+%
\frac{1}{\lambda _{1}^{\left( 2\right) }}\right) ^{3}\left[ \frac{\pi _{1,K}%
}{\xi _{1,2}\left( 1-\pi _{1,K}\right) }+\frac{1}{\mu _{1,2}}+\frac{1-\pi
_{2,K}}{1-\pi _{1,K}}\left( \frac{1}{\lambda _{2}^{\left( 1\right) }}+\frac{1%
}{\lambda _{2}^{\left( 2\right) }}+\frac{1}{\mu _{1,2}}\right) \right] , \\
\pi _{2,K}=\frac{1}{G(NC)}\left( \frac{1}{\lambda _{1}^{\left( 1\right) }}+%
\frac{1}{\lambda _{1}^{\left( 2\right) }}\right) ^{3}\frac{\left( 1-\pi
_{2,K}\right) ^{3}}{\left( 1-\pi _{1,K}\right) ^{3}}\left( \frac{1}{\lambda
_{1}^{\left( 1\right) }}+\frac{1}{\lambda _{1}^{\left( 2\right) }}+\frac{\pi
_{2,K}}{\xi _{2,1}\left( 1-\pi _{1,K}\right) }\right.  \\
\text{ \ \ \ \ \ \ \ }\left. +\frac{1}{\mu _{1,2}}+\frac{1-\pi _{2,K}}{\mu
_{2,1}\left( 1-\pi _{1,K}\right) }\right) ,%
\end{array}%
\right.   \label{eqp-2}
\end{equation}
where the normalization constant $G(NC)$ is given by
\begin{eqnarray}
&&G(NC)  \notag \\
&=&\frac{1}{256\left( \mu _{1,2}\right) ^{4}}+\left( \frac{1}{\lambda
_{1}^{\left( 1\right) }}+\frac{1}{\lambda _{1}^{\left( 2\right) }}\right)
^{3}\frac{\pi _{1,K}}{\xi _{1,2}\left( 1-\pi _{1,K}\right) }+\frac{1}{\mu
_{1,2}}  \notag \\
&&+\frac{1-\pi _{2,K}}{1-\pi _{1,K}}\left( \frac{1}{\lambda _{2}^{\left(
1\right) }}+\frac{1}{\lambda _{2}^{\left( 2\right) }}+\frac{1}{\mu _{2,1}}%
\right)   \notag \\
&&+\left( \frac{1}{\lambda _{1}^{\left( 1\right) }}+\frac{1}{\lambda
_{1}^{\left( 2\right) }}\right) ^{2}\left[ \frac{1}{4\left( \mu
_{1,2}\right) ^{2}}+\frac{\left( 1-\pi _{2,K}\right) ^{2}}{\left( 1-\pi
_{1,K}\right) ^{2}}\left( \frac{1}{\lambda _{2}^{\left( 1\right) }}+\frac{1}{%
\lambda _{2}^{\left( 2\right) }}+\frac{1}{2\mu _{2,1}}\right) ^{2}\right.
\notag \\
&&\left. +\frac{1-\pi _{2,K}}{\mu _{1,2}\left( 1-\pi _{1,K}\right) }\left(
\frac{1}{\lambda _{2}^{\left( 1\right) }}+\frac{1}{\lambda _{2}^{\left(
2\right) }}+\frac{1}{\mu _{2,1}}\right) \right]   \notag \\
&&+\left( \frac{1}{\lambda _{1}^{\left( 1\right) }}+\frac{1}{\lambda
_{1}^{\left( 2\right) }}\right) \left\{ \frac{1}{27\left( \mu _{1,2}\right)
^{3}}+\frac{1-\pi _{2,K}}{4\left( \mu _{1,2}\right) ^{2}\left( 1-\pi
_{1,K}\right) }\left( \frac{1}{\lambda _{2}^{\left( 1\right) }}+\frac{1}{%
\lambda _{2}^{\left( 2\right) }}+\frac{1}{\mu _{1,2}}\right) \right.   \notag
\\
&&\left. +\frac{\left( 1-\pi _{2,K}\right) ^{2}}{\mu _{1,2}\left( 1-\pi
_{1,K}\right) ^{2}}\left( \frac{1}{\lambda _{2}^{\left( 1\right) }}+\frac{1}{%
\lambda _{2}^{\left( 2\right) }}+\frac{1}{2\mu _{1,2}}\right) ^{2}\right.
\notag \\
&&+\frac{\left( 1-\pi _{2,K}\right) ^{3}}{\left( 1-\pi _{1,K}\right) ^{3}}%
\left[ \left( \frac{1}{\lambda _{1}^{\left( 1\right) }}+\frac{1}{\lambda
_{1}^{\left( 2\right) }}\right) ^{3}+\frac{1}{27\left( \mu _{2,1}\right) ^{3}%
}+\frac{1}{4\left( \mu _{2,1}\right) ^{2}}\left( \frac{1}{\lambda
_{2}^{\left( 1\right) }}+\frac{1}{\lambda _{2}^{\left( 2\right) }}\right)
\right.   \notag \\
&&\left. \left. +\frac{1}{\mu _{2,1}}\left( \frac{1}{\lambda _{2}^{\left(
1\right) }}+\frac{1}{\lambda _{2}^{\left( 2\right) }}\right) ^{2}\right]
\right\} +\frac{\left( 1-\pi _{2,K}\right) ^{2}}{4\left( \mu _{1,2}\right)
^{2}\left( 1-\pi _{1,K}\right) ^{2}}\left( \frac{1}{\lambda _{2}^{\left(
1\right) }}+\frac{1}{\lambda _{2}^{\left( 2\right) }}+\frac{1}{2\mu _{2,1}}%
\right) ^{2}  \label{eqp-3} \\
&&+\frac{1-\pi _{2,K}}{3\mu _{1,2}\left( 1-\pi _{1,K}\right) }\left( \frac{1%
}{\lambda _{2}^{\left( 1\right) }}+\frac{1}{\lambda _{2}^{\left( 2\right) }}+%
\frac{1}{\mu _{2,1}}\right)   \notag \\
&&+\frac{\left( 1-\pi _{2,K}\right) ^{3}}{\mu _{1,2}\left( 1-\pi
_{1,K}\right) ^{3}}\left[ \left( \frac{1}{\lambda _{1}^{\left( 1\right) }}+%
\frac{1}{\lambda _{1}^{\left( 2\right) }}\right) ^{3}+\frac{1}{27\left( \mu
_{2,1}\right) ^{3}}+\frac{1}{4\left( \mu _{2,1}\right) ^{2}}\left( \frac{1}{%
\lambda _{2}^{\left( 1\right) }}+\frac{1}{\lambda _{2}^{\left( 2\right) }}%
\right) \right.   \notag \\
&&\left. +\frac{1}{\mu _{2,1}}\left( \frac{1}{\lambda _{2}^{\left( 1\right) }%
}+\frac{1}{\lambda _{2}^{\left( 2\right) }}\right) ^{2}\right]   \notag \\
&&+\frac{\left( 1-\pi _{2,K}\right) ^{4}}{256\left( \mu _{2,1}\right)
^{4}\left( 1-\pi _{1,K}\right) ^{4}}+\frac{\left( \lambda _{2}^{\left(
1\right) }+\lambda _{2}^{\left( 2\right) }\right) \left( 1-\pi _{2,K}\right)
^{4}}{27\lambda _{2}^{\left( 1\right) }\lambda _{2}^{\left( 2\right) }\left(
\mu _{2,1}\right) ^{3}\left( 1-\pi _{1,K}\right) ^{4}}  \notag \\
&&+\frac{\left( \lambda _{2}^{\left( 1\right) }+\lambda _{2}^{\left(
2\right) }\right) ^{2}\left( 1-\pi _{2,K}\right) ^{4}}{4\left( \lambda
_{2}^{\left( 1\right) }\lambda _{2}^{\left( 2\right) }\mu _{2,1}\right)
^{2}\left( 1-\pi _{1,K}\right) ^{4}}  \notag \\
&&+\frac{\left( \lambda _{2}^{\left( 1\right) }+\lambda _{2}^{\left(
2\right) }\right) ^{3}\left( 1-\pi _{2,K}\right) ^{3}\left[ \xi _{2,1}\left(
1-\pi _{2,K}\right) +\pi _{2,K}\mu _{2,1}\right] }{\xi _{2,1}\mu
_{2,1}\left( \lambda _{2}^{\left( 1\right) }\lambda _{2}^{\left( 2\right)
}\right) ^{3}\left( 1-\pi _{1,K}\right) ^{4}}.  \notag
\end{eqnarray}
By using (\ref{eqp-2}) and (\ref{eqp-3}), we can compute the two undetermined constants: $\pi _{1,K}$ and $\pi _{2,K}$.

To this end, let $\lambda _{1}^{\left( 2\right) }=7,\lambda _{2}^{\left( 1\right)
}=5,\lambda _{2}^{\left( 2\right) }=5,\mu _{1,2}=2,\mu _{2,1}=3,\xi
_{1,2}=4,\xi _{2,1}=5$. When $\lambda _{1}^{\left( 1\right) }=5,6,7,8,9$, we obtain the values of $\pi _{1,K}$ and $\pi _{2,K}$ which are listed in Table 1.

From Table 1, it is seen that as $\lambda _{1}^{\left( 1\right) }$ increases, $\pi _{1,K}$ decreases but $\pi _{2,K}$ increases. This result is the same as the actual intuitive situation. When $\lambda _{1}^{\left( 1\right) }$ increases, more bikes are rented from Station 1, so $\pi _{1,K}$ decreases; while when more bikes are rented from Station 1 and are ridden on Road $1\rightarrow 2$, more bikes will be returned to Station 2, so $\pi _{2,K}$ increases.
\begin{table}
\centering
\caption{Numerical results of $\pi _{1,K}$ and $\pi _{2,K}$}
\begin{tabular}{c|c|c|c|c|c|c|c|c|c}
\hline
\thead{$\lambda _{1}^{\left( 1\right) }$}  & \thead{$\lambda _{1}^{\left( 2\right) }$}  & \thead{$\lambda _{2}^{\left( 1\right) }$}  & \thead{$\lambda _{2}^{\left( 2\right) }$} & \thead{$\mu _{1,2} $}& \thead{$\mu _{2,1}$ }& \thead{$\xi_{1,2} $}& \thead{$\xi_{2,1} $}& \thead{$\pi _{1,K} $}& \thead{$\pi _{1,K} $ }\\
\hline
5&7&5&5&2&3&4&5&0.10434&0.14143\\
\hline
6&7&5&5&2&3&4&5&0.08609&0.14502\\
\hline
7&7&5&5&2&3&4&5&0.07609&0.14815\\
\hline
8&7&5&5&2&3&4&5&0.06424&0.14961\\
\hline
9&7&5&5&2&3&4&5&0.05734&0.15116\\
\hline
\end{tabular}
\end{table}

\begin{Rem}
For a large-scale bike sharing system, it is always more difficult and challenging to determine the normalization constant $G(NC)$. Thus it is necessary in the future study to develop some effective algorithms for numerically computing $G(NC)$.
\end{Rem}

\subsection{Performance analysis}

Now, we consider two key performance measures of the bike sharing system in
terms of the steady-state probability $\mathbf{\pi}\left(
\overrightarrow{n}\right)  $ of joint queue lengths at the virtual nodes for $\overrightarrow{n}\in\Omega$.

\textit{(1) The steady-state probability of problematic stations}

In the study of bike sharing systems, it is a key to compute the steady-state probability of problematic stations. For this bike sharing system, the steady-state probability of problematic stations is given by
\begin{align*}
\Im & =\sum\limits_{\text{ }i=1}^{N}P\left\{ n_{i}=0\text{ or }%
n_{i}=K\right\} =\sum\limits_{\text{ }i=1}^{N}\left[ P\left\{
n_{i}=0\right\} +P\left\{ n_{i}=K\right\} \right]  \\
& =\sum\limits_{\text{ }i=1}^{N}\left[ \sum\limits_{\substack{ \text{ }%
\overrightarrow{n}\in \Omega  \\ \&n_{i}=0}}\mathbf{\pi }\left(
\overrightarrow{n}\right) +\sum\limits_{\substack{ \text{ }\overrightarrow{n}%
\in \Omega  \\ \&n_{i}=K}}\mathbf{\pi }\left( \overrightarrow{n}\right) %
\right] .
\end{align*}

\textit{(2) The mean of the steady-state queue length}

The steady-state mean of the number of bikes parked at the $i$th
station is given by%
\begin{equation*}
\mathbf{Q}_{i}=\sum\limits_{\substack{ \text{ }\overrightarrow{n}\in \Omega
\\ \&1\leq n_{i}\leq K}}n_{i}\mathbf{\pi }\left( \overrightarrow{n}\right) ,%
\text{ \ }1\leq i\leq N,
\end{equation*}%
and the steady-state mean of the number of bikes ridden on the Road $ k\rightarrow l$ for $1\leq k\leq N$ and $l\in \Theta _{k}$ is given
by%
\begin{equation*}
\mathbf{Q}_{R_{k\rightarrow l}}=\sum_{r=1,2}\sum_{\substack{ \overrightarrow{%
n}\in \Omega  \\ \&1\leq m_{k,l}^{\left( r\right) }\leq NC}}m_{k,l}^{\left(
r\right) }\mathbf{\pi }\left( \overrightarrow{n}\right) .
\end{equation*}

\begin{Rem}
In the practical bike sharing systems, arrivals of bike users often have some special important behavior and characteristics, such as, time-inhomoge-neity, space-heterogeneity, and arrival burstiness. To express such behavior and characteristics, this paper uses the MAPs to express non-Poisson (and non-renewal ) arrivals of bike users. It is seen that such a MAP-based study is a key to generalize and extend the arrivals of bike users to a more general arrival process in practice, for example, a renewal process, a periodic MAP, a periodic time-inhomogeneous arrival process and so on. In fact, the methodology of this paper may be applied to deal with more general arrivals of bike users. Thus it is very interesting for our future study to analyze the space-heterogeneous or time-inhomogeneous arrivals bike users in the bike sharing systems.
\end{Rem}
\section{Concluding Remarks}

In this paper, we first propose a more general bike sharing system with Markovian arrival processes and under an irreducible path graph.
Then we establish a multiclass closed queueing network by means of some virtual ideas, including, virtual customers, virtual nodes, virtual service times.
Furthermore, we set up the routing matrix, which gives a nonlinear solution to computing the relative arrival rates. Based on this, we give the product-form solution to the steady-state probabilities of joint queue
lengths at the virtual nodes. Finally, we compute the steady-state probability of problematic stations, and also deal with
other interesting performance measures of the bike sharing system. Along these lines, there are a number of interesting
directions for potential future research, for example:

\begin{itemize}
\item Analyzing bike sharing systems with phase type (PH) riding-bike times on
the roads;

\item discussing repositioning bikes by trucks in bike sharing systems with
information technologies;

\item developing effective algorithms for establishing the routing matrix, and for computing the
relative arrival rates;

\item developing effective algorithms for computing the product-form steady-state probabilities of joint queue
lengths at the virtual nodes, and further for calculating the steady-state probability of problematic stations; and

\item applying periodic MAPs, periodic PH distributions, or periodic Markov processes to
study time-inhomogeneous bike sharing systems. This is a very interesting but challenging topic in the future study of bike sharing system.
\end{itemize}
\section*{Acknowledgements}

Q.L. Li was supported by the National Natural Science Foundation of China
under grant No. 71271187 and No. 71471160, and the Fostering Plan of
Innovation Team and Leading Talent in Hebei Universities under grant No. LJRC027.

\end{document}